\input amstex.tex 
\magnification=\magstep{1.5}
\baselineskip 22pt
\documentstyle{amsppt}
\topmatter
\title 
Picard numbers in a family of hyperk\"ahler manifolds \\ 
- A supplement to the article of R. Borcherds, \\ 
L. Katzarkov, T. Pantev, N. I. Shepherd-Barron
\endtitle
\author 
Keiji Oguiso 
\endauthor
\affil 
Department of Mathematical Sciences University of Tokyo
\endaffil
\address
153-8914 Komaba Meguro Tokyo Japan
\endaddress
\email
oguiso\@ms.u-tokyo.ac.jp
\endemail
\dedicatory
Dedicated to Professor Tetsuji Shioda on the occasion of his sixtieth birthday
\enddedicatory
\abstract 
We remark the density of the jumping loci of the Picard number of 
a hyperk\"ahler manifold under small one-dimensional deformation and provide 
some applications for the Mordell-Weil groups of Jacobian K3 surfaces.   
\endabstract
\leftheadtext{Keiji Oguiso}
\rightheadtext{Families of hyperk\"ahler manifolds}
\endtopmatter
\document
\head
{\S 0. Introduction}
\endhead 
Although the N\'eron-Severi groups and the Mordell-Weil groups are described by means of cohomology groups of abelian sheaves, the sheaves are neither coherent nor topological and not much are known about their behaviours under 
deformation. 
\par
In their article ``Families of K3 surfaces'', R. Borcherds, L. Katzarkov, 
T. Pantev and N. I. Shepherd-Barron have found the following phenomenon 
by studying the behaviour of the Picard numbers of K3 surfaces 
under {\it global} deformation: 

\proclaim{Theorem ([BKPS])} Any smooth complete family of minimal K\"ahler 
surfaces 
of Kodaira dimension 0 and constant Picard number is isotrivial. \qed
\endproclaim 

Their proof is a global argument based on the ampleness of the zero 
locus of an automorphic form on a relevant moduli space of K3 surfaces and 
therefore requires the completeness of the base space in essence.
\par
\vskip 4pt

The aim of this note is to remark that their Theorem can be immediately 
generalized to a {\it local one-dimensional} family of hyperk\"ahler 
manifolds if one adopts a different approach (Main Theorem below together 
with Section 1). The idea is extremally simple: In stead of studying behaviours of the Picard groups in family directly, we consider deformation of the 
perpendicular part, i.e. periods. The validity of this reduction is ensured 
by the Lefschetz $(1,1)$-Theorem and the local Torelli Theorem for 
hyperk\"ahler manifolds. We also provide some applications for the Mordell-Weil groups of Jacobian K3 surfaces (Section 2). 
\par
\vskip 4pt

A hyperk\"ahler manifold is by the definition a simply connected, 
compact K\"ahler manifold $F$ with $H^{2, 0}(F) = \Bbb C \omega_{F}$, where 
$\omega_{F}$ is everywhere non-degenerate. According to the Bogomolov 
decomposition Theorem [Be], these are one of the building blocks of manifolds 
with trivial first Chern class. In this terminology, a K3 surface is nothing 
but a hyperk\"ahler manifold of dimension $2$. Due to the fundamental work by 
Bogomolov [Bo] and Beauville [Be], the following results hold for a 
hyperk\"ahler manifold $F$ of any dimension: 
\roster 
\item The Kuranishi space of $F$ is smooth and universal;
\item There exists a primitive integral non-degenerate 
symmetric bilinear form $(*,*)$ on $H^{2}(F, \Bbb Z)$ which induces on 
$H^{2}(F, \Bbb C) = H^{1,1}(F) \oplus \Bbb C\omega_{F} 
\oplus \Bbb C \overline{\omega}_{F}$ the Hodge 
structure of weight two and is of index $(3, B_{2}(F) - 3)$; 
\item The local Torelli Theorem holds for the period map given by the Hodge 
structure on $H^{2}(F, \Bbb Z)$ defined in (2).  
\endroster
Besides original articles, we also refer the readers to [Hu1, Section 1] as 
an excellent survey about these basics. 
\par
\vskip 4pt

In this note, we consider a smooth family of hyperk\"ahler manifolds 
$f : \Cal X \rightarrow \Delta$ over a disk $\Delta$. In this setting, 
the following two statements are equivalent: 
\roster
\item $f$ is trivial as a family, i.e. isomorphic to the product 
$F \times \Delta$ over $\Delta$; 
\item all the fibers of $f$ are isomorphic.
\endroster
This equivalence is a direct consequence of the local Torelli Theorem and the 
universality of the Kuranishi space together with the fact that 
$\pi_{1}(\Delta) = \{1\}$.
\par
\vskip 4pt

We denote by $\rho(F)$ the Picard number of $F$, i.e. the rank of the 
N\'eron-Severi group $NS(F) := \text{Im}(c_{1} : H^{1}(F, \Cal O_{F}^{\times}) 
\rightarrow H^{2}(F, \Bbb Z)) = (\Bbb C \omega_{F})^{\perp} \cap 
H^{2}(F, \Bbb Z)$. Here the last equality is due to the Lefschetz 
$(1,1)$-Theorem. Note also that $0 \leq \rho(F) \leq N := B_{2}(F) - 2$. 
\par
\vskip 4pt

Our main Theorem is as follows: 

\proclaim{Main Theorem} Let $f : \Cal X \rightarrow \Delta$ be a non-trivial 
family of hyperk\"ahler manifolds. Set $M := \text{\rm min}\, \{ \rho(\Cal X_{t}) 
\vert t \in \Delta \}$ and $\Cal S := \{t \in \Delta 
\vert \rho(\Cal X_{t}) > M\}$. Then, $\Cal S$ is a dense countable subset 
of $\Delta$ in the classical topology.
\endproclaim 

The next simple example will illustrate the phenomenon described in the 
main Theorem fairly well:

\example{Example} Let us denote by $E_{t}$ the elliptic curve of 
period $t$. Let $\Delta$ be a small disk in the upper half plane $\Bbb H$. 
Then, one has a family of elliptic curves 
$h : \Cal E \rightarrow \Delta$ with the level two structure 
such that $\Cal E_{t} = E_{t}$. 
Taking a crepant resolution of the quotient of the product 
$g : \Cal E \times E_{\sqrt{-1}} \rightarrow \Delta$ by the inversion, one 
obtains a family of K3 surfaces $f : \Cal X \rightarrow \Delta$ such that 
$\Cal X_{t} = \text{Km}(E_{t} \times E_{\sqrt{-1}})$. This family $f$ satisfies $\rho(\Cal X_{t}) = 20$ for $t \in \Bbb Q(\sqrt{-1})$ and 
$\rho(\Cal X_{t}) = 18$ for 
$t \not\in \Bbb Q(\sqrt{-1})$. In this example, we have 
$\Cal S = \Delta \cap \Bbb Q(\sqrt{-1})$. 
\qed
\endexample 
 
It is an easy fact that $\Cal S$ is at most countable and that the locus of the constant Picard number $\Delta - \Cal S$ is dense and uncountable (and 
therefore ``much bigger'' than $\Cal S$). The essential part of the main 
Theorem is in the converse: {\it existence of enough 
jumping points}. 
\par
\vskip 4pt
As one of applications, we shall solve the following filling up problem of 
possible Picard numbers:
\proclaim{Corollary} Let $F$ be a hyperk\"ahler manifold with $B_{2}(F) 
= N+2$. Then, for each $0 \le j \le N$, there exists a hyperk\"ahler 
manifold $F_{j}$ such that $F$ and $F_{j}$ are deformation equivalent 
and that $\rho(F_{j}) = j$. 
\endproclaim 

\par
\vskip 4pt
Although some form of this kind of density results should be 
known for some experts at least in the case of K3 surfaces (See 
Acknowledgement of [BKPS]), as far as the author knows, there are no 
literatures in which an explicit statement and proof are given. We should 
remark that there are several possible forms of density results and that 
the density of 
the jumping points in a disk implies 
the absence of positive dimensional analytic subsets in the 
Kuranishi space over which the fibers are of constant Picard number. 
Therefore, our statement is stronger than the density of the jumping 
points both in the Kuranishi space and in positive dimensional global base 
spaces. For this remark, it might be worth noticing the following:
\example{Fact} Given a family of K3 surfaces $f : \Cal X \rightarrow \Cal B$, 
we call a point $b \in \Cal B$ a Kummer point if the fiber $\Cal X_{b}$ is 
isomorphic to a Kummer surface. Then: 
\roster 
\item As well-known, the set of Kummer points is dense in the base if the 
family is the 
Kuranishi family of a K3 surface. (See for instance [BPV].) 
\item However, the one dimensional smooth non-trivial family of elliptic 
K3 surfaces $f : \Cal X \rightarrow \Delta_{t}$ ($\vert t \vert << 1$) 
defined by the Weierstrass equation $y^{2} = x^{3} + x + (u^{11} - t)$ 
(over $\Delta_{t}$) has no Kummer points. Indeed, since each fiber 
$\Cal X_{t}$ has non-symplectic automorphism of order $11$ (given by $(x, y, u) \mapsto 
(x, y, \zeta_{11}u)$), the rank of the transcendental lattice must be 
divisible by $\varphi(11) = 10$ and one has then $\rho(\Cal X_{t}) < 16$. 
(See [OZ, Example 2].) \qed
\endroster
\endexample
  
This stronger version is also needed in the proof of the Corollary 
as well as in our application for the 
Mordell-Weil groups of Jacobian K3 surfaces. 
\par
\vskip 4pt 
It would be interesting to ask a similar question for complex tori of 
dimension $\geq 3$. We should also notice that a similar question for 
Calabi-Yau manifolds does not make much sense, because 
$\text{Pic}(X) = H^{2}(X, \Bbb Z)$ for a Calabi-Yau manifold.  
\head Acknowledgement \endhead 
An initial idea of this work was found during the author's stay in the 
Universit\"at 
Essen in 1999 under the financial support by the 
Alexander-Humboldt fundation. Besides [BKPS], a greater part of the idea 
has been brought to the author by Professor Eckart Viehweg through their 
discussion concerning the isotriviality problem [OV]. First of all, the author 
would like to express his deep thanks to Professor Eckart Viehweg for his 
valuable conversation and to the Alexander-Humboldt 
fundation for the financial support. 
Last but not least at all, the author would like to express his 
hearty thanks to Professor Yujiro Kawamata for indicating him several 
improvements, especially, for informing him that the algebraicity assumption 
made in the initial version of this note is unnecessary.    

\head
{\S 1. Proof of the Main Theorem and the Corollary}
\endhead 

Let us choose a marking $\tau : R^{2}f_{*} 
\Bbb Z_{\Cal X} \simeq \Lambda \times \Delta$, where $\Lambda = 
(\Lambda, (*.*))$ is a lattice of signature $(3, N - 1)$ and 
consider the period map 
$$p : \Delta \rightarrow \Cal D := \{ [\omega] \in 
\Bbb P(\Lambda \otimes \Bbb C) \vert (\omega.\omega) = 0, 
(\omega.\overline{\omega}) > 0 \} 
\subset  \Bbb P(\Lambda \otimes \Bbb C) = \Bbb P^{N+1}.$$ 
This map $p$ is defined by $p(t) = \tau_{\Bbb C} ([\omega_{\Cal X_{t}}])$ 
and is known to be holomorphic. We notice that $p$ is not constant by our 
assumption and the local Torelli Theorem. 
\par 
Let us consider all the primitive sublattices $\Lambda_{n}$ ($n \in \Cal N$) 
of $\Lambda$. Put $\Delta_{n} := \{t \in \Delta \vert \tau(NS(\Cal X_{t})) 
= \Lambda_{n} \}$. Then one has a decomposition $\Delta = \sqcup_{n \in \Cal N} \Delta_{n}$. Since $\Cal N$ is countable but 
$\Delta$ is uncountable, there exists an element of $\Cal N$, say $1$, such 
that $\Delta_{1}$ is uncountable. Since 
$p(t) \in \Lambda_{1}^{\perp} \otimes \Bbb C$ for all $t \in \Delta_{1}$ and 
since $p$ is holomorphic, one has then: 
$$p(\Delta) \subset \Cal D' := \{ [\omega] \in 
\Bbb P(\Lambda_{1}^{\perp} \otimes \Bbb C) \vert (\omega.\omega) = 0, 
(\omega.\overline{\omega}) > 0 \} 
\subset  \Bbb P(\Lambda_{1}^{\perp} \otimes \Bbb C) = \Bbb P^{n}.$$ 
Here we regard $\Bbb P(\Lambda_{1}^{\perp} 
\otimes \Bbb C)$ as a linear subspace of $\Bbb P(\Lambda \otimes \Bbb C)$ 
defined by $(\Lambda_{1}.*) = 0$. Set $\Cal S := \Delta - \Delta_{1}$. Then, by the Lefschetz $(1,1)$-Theorem, we also have that:
\roster 
\item $\Lambda_{1} \subset \tau(NS(\Cal X_{t}))$ for all $t \in \Delta$; 
\item $t \in \Cal S$ if and only if there is a vector $v \in \Lambda - 
\Lambda_{1}$ such that $(v, p(t)) = 0$.
\endroster 
Since both $\Lambda - \Lambda_{1}$ and $\{t \in \Delta \vert (v, p(t)) = 0 \}$ 
for each $v \ \in \Lambda - \Lambda_{1}$ are countable, $\Cal S$ is countable 
as well. 
\par
\vskip 4pt 
We shall show the density of $\Cal S$, i.e. the fact that 
$\Cal S \cap U \not= \emptyset$ for any sufficiently small disk $U$. 

\proclaim{Claim 1} $\text{\rm rank}\, \Lambda_{1}^{\perp} \geq 3$. \endproclaim 

\demo{Proof} If $\text{rank} \Lambda_{1}^{\perp} \leq 2$, then 
$\Lambda_{1}^{\perp} \otimes \Bbb R$ is spanned by 
the images of the real and imaginary parts of a holomorphic 2-form. 
This implies that $\Lambda_{1}^{\perp}$ is a positive definite lattice of 
rank $2$ and that $\Cal D'$ consists of two points. However, the period map $p$ is 
then constant, a contradiction. \qed
\enddemo 
Let us choose a holomorphic coordinate $z$ of $U$ centered at $P$. We also 
choose integral  basis of $\Lambda_{1}^{\perp}$ and write $p\vert U$ 
as 
$p(z) = [1 : f_{1}(z) : f_{2}(z) : \cdots : f_{n}(z)]$ with respect to this 
basis. Here we have $n \geq 2$ by Claim 1. We may also assume that $f_{1}(z)$ 
is not constant. In what follows, for each $\overrightarrow{a} = 
(a_{0}, a_{1}, a_{2}, 
\cdots, a_{n}) \in \Bbb R^{n+1} - \{\overrightarrow{0}\}$, we put: 
\par\noindent
$f_{\overrightarrow{a}}(z) := 
a_{0} + a_{1}f_{1}(z) + a_{2}f_{2}(z) + \cdots + a_{n}f_{n}(z)$; 
\par\noindent 
$l_{\overrightarrow{a}} := 
a_{0}x_{0} + a_{1}x_{1} + a_{2}x_{2} + \cdots + a_{n}x_{n}$, where 
$[x_{0} : x_{1} : \cdots : x_{n}]$ is the homogeneous coordinates of 
$\Bbb P^{n}$; and 
\par\noindent
$H_{\overrightarrow{a}} := (l_{\overrightarrow{a}} = 0) \subset \Bbb P^{n}$, 
the hyperplane defined by the linear form $l_{\overrightarrow{a}}$. 
\par 
Let $k$ be an element of $\{2, \cdots, n\}$. Since $\text{dim}_{\Bbb R} 
\Bbb C = 2$, one then finds an element $(r_{0,k}, r_{1,k}, r_{2,k}) 
\in \Bbb R^{3} - \{\overrightarrow{0}\}$ such that 
$r_{0,k}\cdot1 + 
r_{1,k}f_{1}(0) + r_{2,k}f_{k}(0) = 0 - (*)$. 
Put $\overrightarrow{r_{k}} := (r_{0,k}, r_{1,k}, 0, \cdots, 0, r_{2,k}, 0, 
\cdots, 0)$. Here $r_{2,k}$ is located at the same position as $x_{k}$ 
in $[x_{0} : x_{1} : \cdots : x_{n}]$. 
In this notation, the equality $(*)$ is rewritten both as 
$p(0) \in H_{\overrightarrow{r_{k}}}$ and as 
$f_{\overrightarrow{r_{k}}}(0) = 0$. 
\proclaim{Claim 2} $p(U)$ is not contained in 
$\cap_{k = 2}^{n} H_{\overrightarrow{r_{k}}}$. \endproclaim
\demo{Proof} Assuming to the contrary that $p(U) \subset 
\cap_{k = 2}^{n} H_{\overrightarrow{r_{k}}}$, we shall derive a contrdiction. 
Since $f_{1}(z)$ is not constant, we have $r_{2, k} \not = 0$ for each $k$. 
Therefore $\cap_{k = 2}^{n} 
H_{\overrightarrow{r_{k}}}$ is a line $L \simeq \Bbb P^{1}$ defined over 
$\Bbb R$ in $\Bbb P^{n}$. This leads the same contradiction as in Claim 1. 
\qed \enddemo
By Claim 2, there exists $k$ such that $p(U) \not\subset 
H_{\overrightarrow{r_{k}}}$, i.e. $f_{\overrightarrow{r_{k}}}(z) \not\equiv 0$. 
Since $f_{\overrightarrow{r_{k}}}(0) = 0$, we may choose a small circle 
$\gamma \subset U$ around 
$z = 0$ such that $f_{\overrightarrow{r_{k}}}(z)$ has no zeros on $\gamma$. 
Set $K := \min \{\vert f_{\overrightarrow{r_{k}}}(z) \vert \vert z \in 
\gamma\}$ and 
$M := \max \{\vert f_{i}(z) \vert \vert z \in \gamma, i = 0, 1, \cdots, n\}$, 
where we define $f_{0}(z) \equiv 1$. Note that $K > 0$ and $M > 0$. 
Then, by using the triangle inequality, we see that 
$\vert f_{\overrightarrow{r_{k}}}(z) - f_{\overrightarrow{a}}(z) \vert 
< \vert f_{\overrightarrow{r_{k}}}(z) \vert$ on $\gamma$ 
provided that $\vert \overrightarrow{a} - \overrightarrow{r_{k}} \vert < 
KM^{-1}(n+1)^{-1}$. Denote by $V$ the open disk such that $\partial V = 
\gamma$. Then, by the Rouch\'e Theorem, the cardinalities 
of zeros (counted with multiplicities) on $V$ 
are the same for $f_{\overrightarrow{r_{k}}}$ and 
$f_{\overrightarrow{a}}$. In particular, 
$f_{\overrightarrow{a}}$ admits a zero on $V$.
Since $\Bbb Q^{n+1} - \{\overrightarrow{0}\}$ is dense in 
$\Bbb R^{n+1} - \{\overrightarrow{0}\}$, 
one can then find an element $\overrightarrow{q} \in \Bbb Q^{n+1} - 
\{\overrightarrow{0}\}$ such that $f_{\overrightarrow{q}}(z)$ has 
a zero on $V$. Let us denote this zero by $Q \in V (\subset U)$. 
Then, $f_{\overrightarrow{q}}(Q) = 0$ and $p(Q) \in H_{\overrightarrow{q}}$.  
Recall that $\Lambda_{1}^{\perp}$ is primitive in $\Lambda$, $\Lambda$ is 
non-degenerate, and that our homogeneous coordinates 
$[x_{0} : x_{1} : \cdots : 
x_{n}]$ are chosen by means of integral basis of $\Lambda$ and 
the rational linear equations 
$(\Lambda_{1}.*) = 0$. Therefore one can find an element $0 \not = v \in 
\Lambda$ such that $H_{\overrightarrow{q}} = \{x \in \Bbb P^{n} \vert (v.x) 
= 0\}$. Since this $v$ satisfies $(v.p(Q)) = 0$, one has 
$v \in \tau(NS(\Cal X_{Q}))$. On the other hand, since $\overrightarrow{q} 
\not= \overrightarrow{0}$, 
we have $v \not\in \Lambda_{1}$. Hence this $Q$ satisfies $Q \in \Cal S 
\cap U$. Q.E.D.  

\remark{Remark} In general, given a non-constant holomorphic map 
$g : \Delta 
\rightarrow \Bbb P^{n}$, the condition $h(\Delta) \cap H_{\overrightarrow{a}} 
\not = \emptyset$ is not open around $\overrightarrow{a}$ 
if this point satisfies $h(\Delta) \subset H_{\overrightarrow{a}}$. For 
example, the holomorphic map 
$h(z) = [1 : -1 : z]$ satisfies $h(\Delta) \subset H_{(1, 1, 0)}$ and 
$h(\Delta) 
\cap H_{(1-\epsilon, 1, 0)} = \emptyset$ for any $\epsilon \not= 0$. 
\qed
\endremark 

\demo{Proof of Corollary} We shall show by the descending induction on the 
Picard number $\rho := \rho(F)$. Note that $0 \le \rho(F) \le N = 
B_{2}(F) - 2$. Let $u : \Cal U \rightarrow \Cal K$ be the Kuranishi family of 
$F$. This is a germ of 
the universal deformation of $F$ and is of dimension $N$ by 
$H^{1}(T_{F}) \simeq H^{1}(\Omega_{F}^{1})$. 
Threfore, $\Cal K$ is realized as an open neighbourhood of $0 \in 
H^{1}(F, T_{F})$ and is then assumed to be a small polydisk in $\Bbb C^{N}$. 
Then $R^{2}u_{*} \Bbb Z_{\Cal U}$ is a constant system on $\Cal K$. 
Choosing a marking $\tau 
: R^{2}u_{*} \Bbb Z_{\Cal U} \simeq \Lambda \times \Cal K$, 
one can define the period 
map 
$$p :
\Cal K \rightarrow \Cal D := \{ [\omega] \in \Bbb P(\Lambda \otimes 
\Bbb C) \vert (\omega.\omega) = 0, (\omega.\overline{\omega}) > 0 \} 
\subset  \Bbb P(\Lambda \otimes \Bbb C) \simeq \Bbb P^{N+1}.$$ 
This $p$ is a local isomorphism by the local Torelli Theorem and 
by the fact that $\text{dim}\, \Cal K = \text{dim}\, \Cal D = N$. 
Therefore we are allowed to identify $\Cal K$ with a small open set of the
period domain $\Cal D$ (denoted again by $\Cal D$). 
\par
\vskip 4pt
Assume first that $\rho(F) = N$. By $x_{i}$ ($1 \leq i \leq N$), we denote the 
elements in $\Lambda$
corresponding to an integral basis of $\text{NS}(F)$. For each 
$j$ such that $0 \le j < N$, we define the sublocus $\Cal A \subset 
\Cal K = \Cal D$ by the
equations 
$$(x_{1}. \omega) = (x_{2}. \omega) = \cdots =
(x_{j}. \omega) = 0,$$ 
and consider the induced family $\pi :
\Cal S \rightarrow \Cal A$. Here $\Cal A$ is of dimension $N - j > 0$. 
Then, by the construction, we have $\rho(\Cal S_{a}) = j$ 
for generic $a \in \Cal A$ (or more precisely, for any element of the 
complement of the countable union of the hypersurfaces 
$(x.\omega) = 0$, where $x$ runs through the elements of $\Lambda - 
\Bbb Z \langle x_{i} \vert 1 \le i \le j\rangle$), and we are done. 
\par 
\vskip 4pt 
We next assume that $\rho := \rho(F) < N$.  Let us choose an integral 
basis $x_{1}, x_{2}, ..., x_{\rho}$ of $\tau(\text{NS}(F))$. Let us define 
the sublocus $\Cal B \subset \Cal K = \Cal D$ by the
equations 
$$(x_{1}. \omega) = (x_{2}. \omega) = \cdots =
(x_{\rho}. \omega) = 0,$$ 
and consider the induced family $\pi :
\Cal T \rightarrow \Cal B$. Here $\Cal B$ is of dimension $N -
\rho > 0$. By the construction, $\pi$ is not isotrivial and the fibers 
$\Cal T_{b}$
satisfy $\rho(\Cal T_{b}) \geq \rho$. Then by the main Theorem, there is $b
\in \Cal B$ such that $\rho (\Cal T_{b}) > \rho$. Now, by the descending 
induction on $\rho$, we are done. \qed \enddemo
\par
\vskip 4pt 

Exploiting the same idea as in the Example in the Introduction, one also 
immediately obtains the following pretty:

\proclaim{Corollary 2} Set
$\text{GL}^{+}(2, \Bbb Q) := \{M \in \text{GL}(2, \Bbb Q) \vert 
\text{\rm det}\, M > 0 \}$. Let $w = \varphi(z)$ be a holomorphic 
function defined over a neighbourhood of $\tau \in \Bbb H$.
Assume that $\varphi(\tau) \in \Bbb H$. Then, there exists a sequence 
$\{\tau_{k}\}_{k=0}^{\infty} \subset {\Bbb H} - \{\tau\}$ such that 
$\lim_{k \rightarrow \infty} \tau_{k} = \tau$ and that $\tau_{k}$ and 
$\varphi(\tau_{k})$ are congruent for each $k$ under the standard,  
linear fractional action of $\text{GL}^{+}(2, \Bbb Q)$ on $\Bbb H$.  
\endproclaim  

\demo{Proof} Let us denote by $E_{w}$ the elliptic curve of 
period $w$. Choose a small neighbourhood $\varphi(\tau) \in 
\Delta_{2} \subset \Bbb H$. Then one has a family of elliptic curves 
$g : \Cal G \rightarrow \Delta_{2}$ with the level two structure 
such that $\Cal G_{w} = E_{w}$. Since $\varphi$ is holomorphic, we may also 
choose a small neighbourhood $\tau \in \Delta \subset \Bbb H$ such that 
$\varphi(\Delta) \subset \Delta_{2}$ and that there exists a family of 
elliptic curves $h : \Cal H \rightarrow \Delta$ similar to $g$. 
By pulling back $g : \Cal G \rightarrow \Delta_{2}$ by $\varphi$ 
and taking the fiber product, one obtains a family of abelian surfaces 
$a : \Cal A := \Cal H \times_{\Delta} \varphi^{*}\Cal G \rightarrow \Delta$. 
Then, taking a crepant resolution of the quotient of $\Cal A$ by the 
inversion, we obtain a family of Kummer surfaces 
$f : \Cal X \rightarrow \Delta$. By construction, one has $\Cal X_{z} = 
\text{Km}(E_{z} \times E_{\varphi(z)})$. 
\proclaim{Claim} This family 
$f$ is not trivial. \endproclaim 
\demo{Proof} If $f$ is a trivial family, 
then there exists a K3 surface $S$ such that $S \simeq  \text{Km}(E_{z} 
\times E_{\varphi(z)})$ for all $z \in \Delta$. Note that the Kummer surface 
structures on $S$, that is, the isomorphism classes of abelian surfaces $A$ 
such that $S \simeq \text{Km}(A)$, are determined by the choices of $16$ 
disjoint smooth rational curves on $S$. Recall also that there are at most 
countably many smooth rational curves on a K3 surface. Then, there exist 
at most
countably many isomorphism classes of such $A$. Denote all of them by $A_{i}$ 
($i \in \Bbb N$) and set $A_{i} = \Bbb C^{2}/\Lambda_{i}$. For each $A_{i}$, 
the product structures 
on $A_{i}$, i.e. the structures of decompositions $A_{i} = E_{i} \times 
F_{i}$, are also countably many, because subtori of $A_{i}$ are 
determined by the choices of sublattices of $\Lambda_{i}$. Hence, 
there are at most countably many isomorphism classes of pairs 
$(E, F)$ such that $S \simeq \text{Km}(E \times F)$. However, since 
the set of 
the isomorphism classes of $E_{z}$ $(z \in \Delta)$ are uncountable, 
our family $f : \Cal X \rightarrow \Delta$ 
is then non-trivial. \qed 
\enddemo 
Recall by [SM] that $\rho(\Cal X_{z}) = 18$ if 
$E_{z}$ and $E_{\varphi(z)}$ are not isogenous and that 
$\rho(\Cal X_{z}) \geq 19$ if $E_{z}$ and $E_{\varphi(z)}$ are isogenous. 
Then, by the main Theorem, there exists a dense subset $\Cal S \subset \Delta$ 
such that $\rho(\Cal X_{s}) \geq 19$ for $s \in \Cal S$, i.e. that 
$E_{s}$ and $E_{\varphi(s)}$ are isogenous if $s \in \Cal S$. 
Therefore, any sequence in $\Cal S - \{\tau\}$ converging to $\tau$ 
satisfies our requirement. \qed 
\enddemo  

\head 
{\S 2. Appications for Jacobian K3 surfaces} 
\endhead 

A Jacobian K3 surface is an elliptically fibered K3 surface  
$\varphi_{0} : X \rightarrow \Bbb P^{1}$ with section $O$. 
Our interest in this section is the Mordell-Weil group 
$MW(\varphi_{0}) := H^{0}(\Bbb P^{1}, X^{\#})$ of $\varphi_{0}$, i.e. 
the group of sections of $\varphi_{0}$ and its behaviour 
under small deformation. We denote by $r(\varphi_{0})$ the rank of 
$MW(\varphi_{0})$. Note that 
$0 \leq r(\varphi_{0}) \leq 18$ for a Jacobian K3 surface [Sh1]. 
\par 
\vskip 4pt 
By a local one-dimensional family of Jacobian K3 surfaces, we mean a 
commutative diagram 
$$\CD 
\Cal X @> \varphi >> \Cal W @> \Cal O >> \Cal X \\
@V f VV @VV \pi V @VV f V \\
\Delta @> id >> \Delta @> id >> \Delta
\endCD$$ 
such that $\varphi_{t} : \Cal X_{t} \rightarrow \Cal W_{t}$ is a Jacobian 
K3 surface with section $\Cal O_{t}$ for each $t \in \Delta$. We set 
$r(t) := r(\varphi_{t})$. 
\par
\vskip 4pt 
A similar but slightly different jumpinig phenomenon is observed for 
$r(t)$: 

\proclaim{Proposition (2.1)} Let $f = \pi \circ \varphi: 
\Cal X \rightarrow \Cal W \rightarrow \Delta$ be a non-trivial family of 
Jacobian K3 surfaces. Set $\Cal S_{r} := \{t \in \Delta \vert r(t) = r\}$. 
Then, there exists a unique $r_{0}$ such that $\Cal S_{r_{0}}$ is dense and 
uncountable and that 
$\cup_{r > r_{0}} \Cal S_{r}$ is dense and countable. 
\endproclaim 

\remark{Remark} 
\roster 
\item 
The proof below shows that the set $\cup_{r < r_{0}} 
\Cal S_{r}$ is at most countable and has no accumulation point in $\Delta$. 
Moreover, as will be observed in the next Example (2.2), 
there actually exists a case where the set $\cup_{r < r_{0}} 
\Cal S_{r}$ is not empty. Therefore, the behaviour of the Mordell-Weil rank 
in a family is slightly different from that of the Picard number described 
in the main Theorem. 
\item 
As another comparison, we remark that $r(t)$ for a family of rational Jacobian 
surfaces is lower semi-continuous; Therefore the behaviour is quite different 
from the case of Jacobian K3 surfaces described in (2.1). This lower 
semi-continuity is a direct consequence of the stability Theorem and the fact 
that the sections of rational Jacobian surface are $(-1)$-curves. 
\item It would be interesting to study a similar question for a family of 
Jacobian surfaces of Kodaira dimension $1$.  
\qed
\endroster
\endremark  
\demo{Proof} By taking a local trivialization, we may 
assume that $\Cal W = \Bbb P^{1} \times \Delta$. Let $\Cal D \subset 
\Bbb P^{1} 
\times \Delta$ be the discriminant locus of $\varphi$. For the explicit 
description of $\Cal D$, let us consider the Weierstrass model of 
$\varphi : \Cal X 
\rightarrow \Cal W$ (defined by $\Cal O$) and write the equation as 
$y^{2} = x^{3} + a(w, t)x + b(w, t)$, where $w$ is the inhomogeneous 
coordinate of $\Bbb P^{1}$ and $t$ is the coordinate of $\Delta$. Then 
$\Cal D$ is defined by (the reduction of) the equation $4 a(w, t)^{3} + 27 
b(w, t)^{2} = 0 - (*)$. By construction, both $a(w, t)$ and $b(w, t)$ 
are polynomials with respect to $w$. Therefore the restriction map 
$\pi \vert \Cal D : \Cal D 
\rightarrow \Delta; (w, t) \mapsto t$ has at most finitely many such bad 
points $P \in \Cal D$ that $\pi \vert \Cal D$ is not smooth at $P$. Denote by 
$\Cal T \subset \Delta$ the set of the image of these bad points. This is 
then a finite set. In addition, since the type of non-multiple singular 
fibers are uniquely determined by the local monodromy, the singular fibers of 
the fibrations $\varphi_{t} : \Cal X_{t} \rightarrow \Bbb P^{1}$ are 
independent of $t \in \Delta - 
\Cal T$. Write them by $T_{i}$ 
($i = 1, \cdots n$) and denote by $m_{i}$ the number of the irreducible 
components of $T_{i}$. Then, by Shioda's formula [Sh1], we have 
$r(\varphi_{t}) = \rho(\Cal X_{t}) - 2 - \sum_{i = 1}^{n} (m_{i} - 1)$ 
for $t \in \Delta - \Cal T$. Now the result follows from the main Theorem. 
\qed 
\enddemo
 
\example{Example (2.2)} In this example, we shall construct a family 
$f = \pi \circ \varphi: 
\Cal X 
\rightarrow \Cal W \rightarrow \Delta$ of Jacobian K3 surfaces such that 
$\rho(0) = 20$ but $r(0) = 0$, and that $\rho(t) < 20$ but $r(t) > 0$ for 
generic $t$. 
\par 
Let us start from a family of rational 
Jacobian surfaces (with rational double points) $h : \Cal Z \rightarrow 
\Bbb P^{1} \times \Delta_{u} \rightarrow \Delta_{u}$ defined by the 
Weierstrass equation $y^{2} = x^{3} + ux + s^{5}$. Here $u$ is the 
coordinate of $\Delta$ and $s$ is the inhomogeneous coordinate of 
$\Bbb P^{1}$. Then either by the N\'eron algorithm or by a direct 
calculation, one can easily check the following fact: 
$\Cal Z$ is smooth; $\Cal Z_{u}$ ($u \not= 0$) is smooth and 
$\Cal Z_{u} \rightarrow \Bbb P^{1}$ has singular fibers of Type $I_{1}$ over 
$4u^{3} + 27s^{10} = 0$ and of Type $II$ over $s = 
\infty$; and $\Cal Z_{0} \rightarrow \Bbb P^{1}$ has one singular point of 
type $E_{8}$ over $s = 0$ 
and has a singular fiber of Type $II$ over $s = \infty$. 
\par 
Then, by taking an 
appropriate finite covering $\Delta_{v} \rightarrow \Delta_{u}$ and a 
simultaneous resolution of the pull back family, we obtain a family of smooth 
rational Jacobian surfaces $g : \Cal Y \rightarrow \Bbb P^{1} \times 
\Delta_{v} \rightarrow \Delta_{v}$ such that $\Cal Y_{v} \rightarrow 
\Bbb P^{1}$ ($v \not = 0$) has 10 singular fibers of Type $I_{1}$ and one 
singular fiber of Type $II$, and $\Cal Y_{0} \rightarrow \Bbb P^{1}$ has 
one singular fiber of Type $II^{*}$ and one singular fiber of Type $II$. 
By Shioda's formula, one has then $r(v) = 8$ for $v \not= 0$ and $r(0) = 0$. 
\par
Let us choose large number $M$ such that the divisor $s = M$ on 
$\Bbb P^{1} \times \Delta_{v}$ does not meet the discriminant locus. 
This is possible by the description above. Let us take the double covering 
$\Bbb P^{1} \times \Delta_{t} \rightarrow \Bbb P^{1} \times \Delta_{v}$ 
ramified over $s = M$ and $s = \infty$ and consider the relatively 
minimal model $f : \Cal X \rightarrow \Bbb P^{1} \times \Delta_{t} 
\rightarrow \Delta_{t}$ of the pull back family 
$f' : \Cal X' \rightarrow \Bbb P^{1} \times \Delta_{t} \rightarrow \Delta_{t}$. Note that $\Cal X'$ is equi-singular along the preimage of 
the cuspidal points of fibers $(\Cal Y_{v})_{\infty}$ 
of $\Cal Y_{v} \rightarrow \Bbb P^{1}$. Therefore, by the monodromy 
calculation, one finds that $f : \Cal X \rightarrow \Bbb P^{1} \times 
\Delta_{t} 
\rightarrow \Delta_{t}$ is a smooth family of Jacobian K3 surfaces such that 
$\Cal X_{0} \rightarrow \Bbb P^{1}$ has two singuler fibers of Type $II^{*}$ 
and one singular fiber of Type $IV$; and $\Cal X_{t} \rightarrow \Bbb P^{1}$ 
($t \not= 0$) has 20 singular fibers of Type $I_{1}$ and one singular 
fiber of Type $IV$. Note also that $r(t) \geq r(v) = 8$ for $t \not= 0$. 
On the 
other hand, again by Shioda's formula, one has $20 \geq \rho(t=0) = 2 + 
r(t=0) + (9 -1) + (9 -1) + (3 -1)$. Therefore $\rho(0) = 20$ and 
$r(0) = 0$ for the 
central fiber of $f$. Moreover, this family $\Cal X \rightarrow \Delta$ is not 
trivial as a family of K3 surfaces. (Indeed, otherwise, we have $\Cal X 
\simeq \Cal X_{0} \times \Delta_{t}$. Since $\text{Pic}(\Cal X_{0})$ is 
discrete, $\Cal X_{0}$ does not admit a family of elliptic 
fibrations varying continuously. Then, our family 
$\Cal X \rightarrow \Delta$ must be also trivial 
as a family of elliptic fiber spaces. However, this contradicts the fact 
that the types of singular fibers of $\Cal X_{0}$ and $\Cal X_{t}$ are 
different.) Hence, by the main Theorem, we have $\rho(t) < 20$ 
for generic $t$. 
\qed 
\endexample
\remark{Remark} This example also shows that there is a case where the 
behaviour of $r(t)$ is not honestly accompanied with that of $\rho(t)$. \qed
\endremark
Next we apply our main Theorem to study the structure of 
the Mordell-Weil lattices of Jacobian K3 surfaces. 
Here the Mordell-Weil group $MW(\varphi)$ with Shioda's positive definite, 
symmetric bilinear form $\langle *, *\rangle$ is called 
the Mordell-Weil lattice [Sh2]. 
This lattice structure on $MW(\varphi)$ made the study of Mordell-Weil groups 
extremally rich [Sh3]. By the narrow Mordell-Weil 
lattice $MW^{0}(\varphi)$ we mean the sublattice of $MW(\varphi)$ of finite 
index consisting of the sections which pass through the identity component of 
each fiber [Sh2]. Contrary to the case of rational Jacobian 
surfaces, the isomorphism classes of both $MW(\varphi)$ and $MW^{0}(\varphi)$ 
for Jacobian K3 surfaces are no more finite ([OS], [Ni]) and the whole 
pictures of them does not seem so clear even now. Our interest here is 
to clarify certain relationships among all of the Mordell-Weil lattices 
of Jacobian K3 surfaces:

\proclaim{Theorem (2.3)} For any given Jacobian K3 surface $\varphi : J 
\rightarrow \Bbb P^{1}$ of rank $r := r(\varphi)$, there exists a sequence 
$\{\varphi_{m} : J_{m} \rightarrow \Bbb P^{1}\}_{m = r}^{18}$ of Jacobian K3 
surfaces such that 
\roster
\item 
$\varphi_{r} : J_{r} \rightarrow \Bbb P^{1}$ is 
the original $\varphi : J \rightarrow \Bbb P^{1}$;
\item 
$r(\varphi_{m}) = m$ for each $m$; and 
\item 
there exists a sequence of isometric embeddings:
\endroster
$$MW^{0}(\varphi) (= MW^{0}(\varphi_{r})) \subset MW^{0}(\varphi_{r+1}) 
\subset \cdots \subset MW^{0}(\varphi_{17}) \subset MW^{0}(\varphi_{18}).$$   
In particular, the narrow Mordell-Weil lattice of a Jacobian K3 surface is 
embedded into the Mordell-Weil lattice of some Jacobian K3 
surface of rank $18$. Conversely, for every sublattice $M$ of the (narrow) 
Mordell-Weil lattice of a Jacobian K3 surface of rank $18$, there exists 
a Jacobian K3 surface whose narrow Mordell-Weil lattice contains 
$M$ as a sublattice of finite index. Moreover, for each given $M$ 
there are at most finitely many isomorphism classes of the Mordell-Weil 
lattices of Jacobian K3 surfaces which contains $M$ as a sublattice of 
finite index.  
\endproclaim 

\demo{Proof} First we shall show the existence of a sequence 
in the statement. We may assume that $r \leq 17$. 
Let us consider the Kuranishi family 
$k : (J \subset \Cal U) \rightarrow (0 \in \Cal K)$ of $J$. This is a germ of 
the universal deformation of $J$ and is known to be smooth of dimension 20. 
Threfore, $\Cal K$ is realized as an open neighbourhood of $0 \in 
H^{1}(J, T_{J})$ and is then assumed to be a small polydisk in $\Bbb C^{20}$. 
Then $R^{2}k_{*} \Bbb Z_{\Cal U}$ is a constant system on $\Cal K$. 
Choosing a marking $\tau 
: R^{2}k_{*} \Bbb Z_{\Cal U} \simeq \Lambda \times \Cal K$, where 
$\Lambda = \Lambda_{\text{K3}}$, let us consider as before 
the period map 
$$p : \Cal K \rightarrow \Cal D = \{ [\omega] \in \Bbb P(\Lambda \otimes 
\Bbb C) \vert (\omega.\omega) = 0, (\omega.\overline{\omega}) > 0 \} 
\subset  \Bbb P(\Lambda \otimes \Bbb C) \simeq \Bbb P^{21}.$$ 
Since $p$ is a local isomorphism for the same reason as before, 
we may identify $\Cal K$ with an open neighbourhood $\Cal U \subset \Cal D$ 
of $p(0)$ by $p$. Since our argument is completely local, 
by abuse of notation, we 
write this $\Cal U$ again by $\Cal D$ and identify therefore $\Cal K = \Cal D$ 
by $p$. 
\par 
Let us write a general fiber of $\varphi : J \rightarrow 
\Bbb P^{1}$ by $E$ and choose an integral basis $S_{i}$ ($i = 1, \cdots, r$) 
of $MW^{0}(\varphi)$, where $r := r(\varphi)$. Then $S_{i}$ and the zero 
section $O$ are all non-singular rational curves and $E$ is an elliptic curve 
such that $(S_{i}.E) = (O.E) = 1$. By the definition of the Mordell-Weil 
lattice 
$(MW(\varphi), \langle *, * \rangle)$, we have a (minus sign of) 
isometric injective homomorphism  $\iota : MW^{0}(\varphi) \hookrightarrow 
NS(J)$ given by $S \mapsto S - O + ((O^{2}) - (S.O))E$ [Sh2]. Then, 
$\langle S_{i}, S_{i} 
\rangle = 4 + 2(S_{i}.O)$. Note also that $E$, $O$, $S_{i}$ are linearly 
independent in $NS(J)$ [Sh1]. Let us consider $(r+2)$ elements in 
$\Lambda$ given by $e := \tau([E])$, $o := \tau([O])$ and 
$s_{i} := \tau([S_{i}])$. Then, these are also linearly indepedent in 
$\Lambda$. 
\par
\vskip 4pt 
Consider the subset $\Cal L$ of $\Cal D$ defined by 
$(e.*) = (o.*) = (s_{i}.*) = 0$. This is a smooth analytic subset of 
$\Cal D$ of dimension $20 - (r + 2) > 0$ and contains $p(0)$. Through the 
identification made 
above, we may regard $0 \in \Cal L \subset \Cal K$. Then we may speak of the 
family $\tilde{j} : \tilde{\Cal J} \rightarrow \Cal L$ obtained as the 
restriction of 
$k : \Cal U \rightarrow \Cal K$ to $\Cal L$. Then by [Ko, Theorem 14] or 
by [Hu1, Section 1 (1.14)], one finds that 
$\Cal L$ is the locus on which the invertible sheaves $\Cal O_{J}(E)$, 
$\Cal O_{J}(O)$, $\Cal O_{J}(S_{i})$ on $J$ lift to invertible sheaves 
$\Cal E$, $\Cal O$ and $\Cal S_{i}$ on the whole space $\tilde{\Cal J}$. 
Since $20 - (r+2) \geq 1$ by $r \leq 17$, one can take a 
sufficiently small disk $0 \in \Delta \subset \Cal L$ and obtains  
the induced family $j : \Cal J \rightarrow \Delta$. We denote the restrictions 
of $\Cal E$, $\Cal O$ and $\Cal S_{i}$ on $\Cal J$ by the same letter. We also 
shrink $\Delta$ freely whenever it is convenient. 
Note that $\chi(\Cal O_{J}(S_{i})) = 1$, $h^{0}(\Cal O_{J}(S_{i})) = 1$ 
and $h^{q}(\Cal O_{J}(S_{i})) = 0$ for $q > 0$, because $S_{i}$ is a smooth 
rational curve on a K3 surface. Then by applying the upper 
semi-continuity of coherent sheaves and by the Theorem of cohomology, we 
see that $j_{*}\Cal S_{i}$ are invertible sheaves which satisfy the base 
change property. Then $(j_{*}\Cal S_{i}) \otimes \Bbb C(0) 
\simeq H^{0}(\Cal O_{J}(S_{i}))$. Therefore, by Nakayama's Lemma, 
all of $C_{i}$ lift not only as invertible sheaves but also as effective 
divisors on $\Cal J$. By abuse of notation, we denote these divisors again 
by $\Cal S_{i}$. Since the smoothness is an open condition for a proper 
morphism, 
$\pi \vert \Cal S_{i} : \Cal S_{i} \rightarrow \Delta$ is also smooth. 
Combining this with the fact that small 
deformation of $\Bbb P^{1}$ is again $\Bbb P^{1}$, we see that $S_{i, t} := 
\Cal S_{i} \vert \Cal J_{t}$ is again a smooth rational curve on $\Cal J_{t}$ 
for all $t \in \Delta$. The same holds for $\Cal O_{t} := \Cal O \vert 
\Cal J_{t}$. Note that $\chi(\Cal O_{J}(E)) = 2$, $h^{0}(\Cal O_{J}(E)) = 2$ 
and $h^{q}(\Cal O_{J}(E)) = 0$ for $q > 0$, because $E$ is an elliptic curve 
on a K3 surface. Then, $h^{q}(\Cal E \vert \Cal J_{t}) = 0$ and 
$h^{0}(\Cal E \vert \Cal J_{t}) = 2$. Therefore, $j_{*} \Cal E$ is a locally 
free sheaf of rank 2 which satisfies the base change property. In 
particular, $j^{*}j_{*}\Cal E \vert J = H^{0}(\Cal O_{J}(E))$. Since 
$\Cal O_{J}(E)$ is globally generated, we see again by Nakayama's Lemma 
that the natural map $j^{*}j_{*}\Cal E \rightarrow \Cal E$ is also surjective. 
Therefore we may 
take a morphism $\Phi : \Cal J \rightarrow \Cal W$ over $\Delta$ associated 
to this surjection. Then, by the base change property, we find that the 
restriction $\Phi_{t} : \Cal J_{t} \rightarrow 
\Cal W_{t}$ coincides with the morphism given by the surjection 
$H^{0}(\Cal E \vert \Cal J_{t}) \otimes \Cal O_{\Cal J_{t}} \rightarrow 
\Cal E \vert \Cal J_{t}$. This is an elliptic fibration by 
$h^{0}(\Cal E \vert \Cal J_{t}) = 2$ and by the adjunction formula on a 
K3 surface. Therefore, the factorization $\Phi : \Cal J \rightarrow \Cal W$ 
makes $j : \Cal J \rightarrow \Delta$ a family of elliptic surfaces over 
$\Delta$. By the invariance of the intersection number, we have 
$(\Cal S_{i, t}. \Cal E_{t}) = (S_{i}.E) = 1$. Therefore, $S_{i,t}$ is also 
a section of $\Phi_{t}$. The same holds for $\Cal O_{t}$. Therefore 
$\Phi : \Cal J \rightarrow \Cal W$ makes $j : \Cal J \rightarrow \Delta$ 
a family of Jacobian K3 surfaces with zero section $\Cal O$. Moreover, 
by passing to 
the Weierstrass family over $\Delta$ given by $\Cal O$ and using the 
characterization of $MW^{0}(\varphi)$ that $S \in MW(\varphi)$ is in 
$MW^{0}(\varphi)$ if and only if $S$ does not meet the singular points of the 
Weierstrass model, one finds that $\Cal S_{i, t}$ are all in 
$MW^{0}(\Phi_{t})$. In addition, the 
intersection matrix of $\Cal E_{t}$, $\Cal O_{t}$, $\Cal S_{i, t}$ are the 
same as the one for $E$, $O$, $S_{i}$ in $\Lambda$ and is then hyperbolic. 
Therefore, $\Cal E_{t}$, $\Cal O_{t}$, $\Cal S_{i, t}$ are also linearly 
independent in $H^{2}(\Cal J_{t}, \Bbb Z)$. Hence so are in $NS(\Cal J_{t})$. 
Thus by the 
injection $MW^{0}(\Phi_{t})  \hookrightarrow NS(\Cal J_{t})$ quoted above, we 
see that $S_{i,t}$ are also linearly independent in $MW^{0}(\Phi_{t})$. In 
particular, $r(\Phi_{t}) \geq r$ for all $t \in \Delta$. Since the base space 
$\Delta$ is chosen in the Kuranishi space, our family $j : \Cal J \rightarrow 
\Delta$ is not trivial. Therefore, by Proposition (2.1), there exists $t_{0} 
\in 
\Delta$ such that $r(t_{0}) > r$. By the invariance intersection and 
by the relation between $\langle *, * \rangle$ and $(*, *)$ quoted above, 
we see that the map $a : MW^{0}(\varphi) 
\rightarrow MW^{0}(\Phi_{t_{0}})$ given by $S_{i} \mapsto \Cal S_{i, t_{0}}$ 
is then an isometric injection. 
\par 
If $r(t_{0}) = r + 1$, then we may define 
$\varphi_{r+1} : J_{r+1} \rightarrow \Bbb P^{1}$ to be this Jacobian K3 
surface $\Phi_{t_{0}} : \Cal J_{t_{0}} 
\rightarrow \Bbb P^{1}$. 
\par 
Let us treat the case where 
$r(t_{0}) \geq r + 2$. Since $18 \geq r(t_{0})$, we have $16 \geq r$. 
For simplicity, 
we abbreviate $\Phi_{t_{0}} : \Cal J_{t_{0}} 
\rightarrow \Bbb P^{1}$ and $r(t_{0})$ by 
$\varphi' : J' \rightarrow \Bbb P^{1}$ and $r'$ respectively. We denote the 
image of the basis $S_{i}$ ($1 \leq i \leq r$) of $MW^{0}(\varphi)$ 
in $MW^{0}(\varphi')$ by the same letters $S_{i}$ and take 
$T_{j} \in MW^{0}(\varphi')$ $j = r+1, r+2, \cdots, r'$ such that 
$S_{i}$ and $T_{j}$ form a basis of $MW^{0}(\varphi') \otimes 
\Bbb Q$ over $\Bbb Q$. (Here note that our embedding might not be primitive 
so that we can not prolong $S_{i}$ to integral basis of 
$MW^{0}(\varphi')$ in general.) Let us consider the Kuranishi space 
$\Cal K'$ of $J'$ and take the subspace $\Cal L' \subset \Cal K'$ 
defined by the fiber  
class $E'$ of $\varphi'$, the zero section $O$, all of $S_{i}$, and 
$T_{r+1}$. Denote by $j' : \Cal J' \rightarrow \Cal L'$ the family 
induced by the Kuranishi family as before. Then, 
$\text{dim}\Cal L' = 20 - (2 + r + 1) > 0$, because $E'$, $O$, 
$S_{i}$ and $T_{r+1}$ are linearly independent in $H^{2}(J', \Bbb Z)$ 
and $r \leq 16$. 
In addition, considering $\Cal L'$ as a 
subspace in the period domain under the identification made as before, and 
applying the same argument as in the proof of the main Theorem, 
one finds that the 
N\'eron-Severi group of $\Cal J_{t}'$ for $t$ being generic in $\Cal L'$ is 
isomorphic to the primitive 
closure of $\Bbb Z \langle E_{18}, O, S_{i}, T_{r+1} \rangle$ in 
$H^{2}(J', \Bbb Z)$. In particular, $\rho(\Cal J_{t}') = r + 3$ 
for generic $t$. Moreover, 
by the same argument as above, one makes this family a 
family of Jacobian K3 surfaces 
$\Cal J' \overset \Phi' \to \longrightarrow \Cal W' \rightarrow \Cal L'$ 
such that each fiber $\Phi_{t}' : \Cal J_{t}' \rightarrow \Cal W_{t}'$ 
satisfies that $MW^{0}(\varphi) 
\subset MW^{0}(\Phi_{t}')$ and that
$r(\Phi_{t}') \geq r + 1$. On the other hand, one has $r(\Phi_{t}') \leq 
r + 1$ for generic $t$ by Shioda's formula and by $\rho(\Cal J_{t}') = r + 3$. 
Then we have $r(\Phi_{t}') = r + 1$ and may define $\varphi_{r+1} : J_{r +1} 
\rightarrow \Bbb P^{1}$ to be $\Phi_{t}' : \Cal J_{t}' \rightarrow 
\Cal W_{t}'$ for generic $t$. 
The first statement now follows from induction on $(18 - r)$. 
\par
Next we shall show the middle statement. Let $\phi' : S' 
\rightarrow \Bbb P^{1}$ be a Jacobian K3 
surface such that $\text{rank}(\phi') = 18$ and $M$ a sublattice of 
$MW(\phi')$. Then, by taking a generic point of the locus of the Kuranishi 
space defined by the basis of $M$, zero section of 
$\phi'$ and general fiber of $\phi'$, one gets a Jacobian K3 
surface 
$\phi : S \rightarrow \Bbb P^{1}$ such that $M \subset MW^{0}(\phi)$ 
and $r(\phi) = r$. 
\par 
Finally, we check the last assertion. 
Assume that $M \subset MW^{0}(\phi)$ and is of finite index. 
Since the pairing $\langle * , * \rangle$ is integral valued on 
$MW^{0}(\phi)$ [Sh2], one has $M \subset MW^{0}(\phi) \subset M^{*}$. 
Since $M \subset M^{*}$ is of finite index, the possibilities of 
$MW^{0}(\phi)$ is then only finitely many. By [Sh2], one has also 
$MW^{0}(\phi) \subset MW(\phi)/\text{(torsion)} \subset 
MW^{0}(\phi)^{*}$. Therefore each $MW^{0}(\phi)$ recovers 
$MW(\phi)/\text{(torsion)}$ up to finitely many ambiguities. 
Now it is suffices to check the boundedness of the torsion subgroups of 
Jacobian K3 surfaces. If the $j$-invariant is not constant, 
the result 
follows from the classification due to Cox [Co]. Let us consider the case 
where the $j$-invariant is constant. Note that a Jacobian 
K3 surface always admits at least one singular fiber, because 
its topological Euler number is positive. Therefore, by the classification of 
the singular fibers whose $j$-values are not $\infty$ and 
by the general fact that 
the specialization map $MW(\phi)_{\text{torsion}} \rightarrow 
(\phi^{-1}(t))_{\text{reg}}$ is 
injective, one can easily see that the possible torsion groups are at most 
$0$, $\Bbb Z/2$, $\Bbb Z/3$, $\Bbb Z/4$ or $(\Bbb Z/2)^{\oplus 2}$. 
Now we are done.  
\qed 
\enddemo

This Theorem coarsely reduces the study of $MW(\varphi)$ to those of 
the maximal rank $18$. For further study, it might be worthwhile noticing that 
a Jacobian K3 surface with maximal Mordell-Weil rank 18 is necessarily 
``singular'' in the sense of Shioda [SI] and that Nishiyama [Ni] has already 
constructed an infinite series of examples of such Jacobian K3 surfaces. 

\enddocument